\documentclass[12pt,a4paper,draft]{article}
\usepackage[english,russian]{babel}
\usepackage{amssymb, amsmath, latexsym, amsfonts, amsthm}
\begin{document}
\renewcommand{\refname}{References}
\newtheorem{theorem}{Theorem}
\newtheorem{lemma}{Lemma}
\newtheorem{corollary}{Corollary}
\newtheorem{proposition}{Proposition}
\begin{center}
On a space of entire functions rapidly decreasing on ${\mathbb R}^n$ and its Fourier transformation
\end{center}
\begin{center}
I. Kh. Musin\footnote {E-mail: musin@matem.anrb.ru}, M. I. Musin\footnote {E-mail: marat402@gmail.com}
\end{center}
\begin{center}
Institute of Mathematics with Computer Centre of Ufa Scientific Centre of Russian Academy of Sciences, 

Chernyshevsky str., 112, Ufa, 450077, Russia;

Bashkirian state University, Department of Mathematics and Computer technologies, 

Z. Validy str., 32, Ufa, 450000, Russia
\end{center}

\vspace {0.3cm}

\renewcommand{\abstractname}{}
\begin{abstract}
{\sc Abstract}. A space of entire functions of several complex variables rapidly decreasing on ${\mathbb R}^n$ and such that their growth along $i{\mathbb R}^n$ is majorized with a help of a family of weight functions is considered in the paper. For this space an equivalent description in terms of estimates on all partial derivatives of functions on ${\mathbb R}^n$ and Paley-Wiener type theorem are obtained. 

\vspace {0.3cm}






\vspace {0.3cm}

\end{abstract}

\vspace {1cm}
 
\section{Introduction}

{\bf 1.1. On the problem}. 
Let $\varPhi =\{\varphi_m\}_{m=1}^{\infty}$ be a family of nondecreasing continuous functions 
$\varphi_m: [0, \infty) \to {\mathbb R}$
such that for each $m \in {\mathbb N}$:

$i_1)$. $\displaystyle \lim_{x \to + \infty} \frac {\varphi_m(x)}{x}= + \infty$;

$i_2)$. for each $A > 0$ there exists a constant $C(m, A) > 0$ such that 
$$
\varphi_m(x) + A \ln (1 + x) \le \varphi_{m+1}(x) + C(m, A), \ x \ge 0.
$$

Let $H({\mathbb C}^n)$ be the space of entire functions on ${\mathbb C}^n$,  
$\Vert u \Vert$ be the Euclidean norm of $u \in {\mathbb R}^n ({\mathbb C}^n)$. 
For each $\nu \in {\mathbb N}$ and $k \in {\mathbb Z}_+$ define the Banach space
$$
E_k(\varphi_{\nu}) = \{f \in H({\mathbb C}^n): p_{\nu, k}(f) = \sup_{z \in {\mathbb C}^n} 
\frac 
{\vert f(z)\vert (1 + \Vert z \Vert)^k}
{e^{\varphi_{\nu} (\Vert Im \ z \Vert)}} < \infty \}.
$$
Let $E(\varphi_{\nu})= \bigcap \limits_{k=0}^{\infty} E_k(\varphi_{\nu})$, $E(\varPhi)= \bigcup \limits_{\nu=1}^{\infty} E(\varphi_{\nu})$. 
With usual operations of addition and multiplication by complex numbers 
$E(\varphi_{\nu})$ and $E(\varPhi)$ are linear spaces. 
Since $p_{\nu, k}(f) \le p_{\nu, k+1}(f)$ for $f \in E_{k+1}(\varphi_{\nu})$, then $E_{k+1}(\varphi_{\nu})$ is continuously embedded in $E_k(\varphi_{\nu})$. 
Endow $E(\varphi_{\nu})$ with a projective limit topology of spaces $E_k(\varphi_{\nu})$. 
Note that if $f \in E(\varphi_{\nu})$ then 
$p_{\nu+1, k}(f) \le e^{C(\nu, 1)} p_{\nu, k}(f)$ for each $k \in {\mathbb Z}_+$. 
This means that $E(\varphi_{\nu})$ is continuously embedded in $E(\varphi_{\nu + 1})$ for each $\nu \in {\mathbb N}$. 
Supply $E(\varPhi)$ with a topology of an inductive limit of spaces $E(\varphi_{\nu})$.

Note that if functions $\varphi_{\nu}$ are defined by the formula 
$\varphi_{\nu}(x) = \Omega(\nu x)$ ($\nu \in {\mathbb N}$) where $\Omega$ be differentiable functions on $[0, \infty)$ such that $M(0)= \Omega(0) = M'(0) = \Omega'(0) = 0$ and their derivatives are continuous, increasing and tending to infinity then $E(\varPhi)$ is Gelfand-Shilov space $W^{\Omega}$ \cite {Gur1}, \cite {Gur2}, \cite {GS1}, \cite {GS2}, \cite{C-C-K 1}, \cite   {C-C-K 2}.

In the paper we describe the space $E(\varPhi)$ in terms of estimates 
on partial derivatives of functions on ${\mathbb R}^n$ and study Fourier transformation of functions of $E(\varPhi)$ under rather weak additional conditions on $\varPhi$. 

{\bf 1.2. Some notations 	and definitions}.  
For $u=(u_1, \ldots , u_n) \in {\mathbb R}^n \ ({\mathbb C}^n)$, $v=(v_1, \ldots , v_n) \in {\mathbb R}^n \ ({\mathbb C}^n)$ \ 
$\langle u, v \rangle  = u_1 v_1 + \cdots + u_n v_n$. 

For $\alpha = (\alpha_1, \ldots , \alpha_n) \in {\mathbb Z}_+^n$, 
$x =(x_1, \ldots , x_n) \in {\mathbb R}^n$, 
$z =(z_1, \ldots , z_n) \in {\mathbb C}^n$ \ $\vert \alpha \vert = \alpha_1 + \ldots  + \alpha_n$, 
$x^{\alpha} = x_1^{\alpha_1} \cdots x_n^{\alpha_n}$, $z^{\alpha} = z_1^{\alpha_1} \cdots z_n^{\alpha_n}$, 
$D^{\alpha}=
\frac {{\partial}^{\vert \alpha \vert}}{\partial x_1^{\alpha_1} \cdots \partial x_n^{\alpha_n}}$ .

For multi-indices
$\alpha = ({\alpha}_1, \ldots , {\alpha}_n),
\beta = ({\beta}_1, \ldots , {\beta}_n) \in {\mathbb Z_+^n}$ the notation
$\beta \le \alpha $ indicates that
${\beta}_j \le {\alpha}_j$ ($j = 1, 2, \ldots , n$). 

For multi-indices
$\alpha = ({\alpha}_1, \ldots , {\alpha}_n),
\beta = ({\beta}_1, \ldots , {\beta}_n) \in {\mathbb Z_+^n}$ such that
$\beta \le \alpha $ let $C_{\alpha}^{\beta}= \prod \limits_{j=1}^{n} C_{\alpha_j}^{\beta_j}$, where $C_{\alpha_j}^{\beta_j}$ are the combinatorial numbers. 

By $s_n(1)$ denote the surface area of the unit sphere in ${\mathbb R}^n$.

For a function $u: [0, \infty) \to {\mathbb R}$ let
$u[e](x): = u(e^x), \ x \ge 0$.

For brevity we put $\psi_m:=\varphi_m [e]$ ($m \in {\mathbb N}$).

Denote by ${\mathcal B}$ the set of all real-valued functions $g \in C[0, \infty)$ such that 
$\displaystyle \lim_{x \to + \infty} \frac {g(x)}{x}=~+\infty$. 

Let 
$V =\{h \in {\mathcal B}: \text {$h$ is convex on $[0, \infty)$}\}$,
${\mathcal V} =\{h \in V: \text {$h$ is increasing on $[0, \infty)$ with $h(0)=0$} \}$.

For $g \in V$ let 
$
V_g =\{h \in {\mathcal V}: \text {$h$ coincides with 
$g$ on $[d_h, \infty)$},
$
where $d_h$ is some positive number depending on $h$\}.


For $g \in {\mathcal B}$ let $g^*$ be its the Young conjugate:
$g^*(x) = \displaystyle \sup \limits_{y \ge 0}(xy - g(y)), \ x \ge 0$. 
Recall that by the inversion formula for Young transform \cite {R-V} if $ g \in {\mathcal V}$ 
then $(g^*)^*=g$.

{\bf 1.3. Main results}. Let $\Psi^*=\{\psi_{\nu}^*\}_{\nu=1}^{\infty}$. For each $\nu \in {\mathbb N}$ and $m \in {\mathbb Z}_+$ let 
$$
{\mathcal E}_m(\psi_{\nu}^*) =\{f \in  C^{\infty}({\mathbb R}^n): 
{\cal R}_{m, \nu}(f)= \sup_{x \in {\mathbb R}^n, \alpha \in {\mathbb Z}_+^n} 
\frac {(1+ \Vert x \Vert)^m \vert (D^{\alpha}f)(x) \vert}{\alpha! e^{-\psi_{\nu}^*(\vert \alpha \vert)}} < \infty \}.
$$
Let ${\mathcal E}(\psi_{\nu}^*) = \bigcap \limits_{m=0}^{\infty}{\mathcal E}_m(\psi_{\nu}^*)$, 
${\mathcal E}(\Psi^*) = \bigcup \limits_{\nu=1}^{\infty}{\mathcal E}(\psi_{\nu}^*)$.

The first two theorems (proved in section 3) are aimed to characterize functions of the space $E(\varPhi)$ 
in terms of estimates of their partial derivatives on ${\mathbb R}^n$.   

\begin{theorem} 
Let the family $\varPhi$ satisfies the additional condition

$i_3)$. for each $m \in {\mathbb N}$ there is a constant $a_m > 0$ such that 
$$
\varphi_m(2x) \le \varphi_{m+1}(x) + a_m, \ x \ge 0.
$$

Then for each $f \in  E(\varPhi)$ its restriction $f_{|{\mathbb R}^n}$ 
on ${\mathbb R}^n$ belongs to ${\mathcal E}(\Psi^*)$.
\end{theorem}

\begin{theorem} 
Let the family $\varPhi$ satisfies the additional condition

$i_4)$. for each $m \in {\mathbb N}$ there exist constants $\sigma_m > 1$ and $\gamma_m > 0$ such that
$$
\varphi_m(\sigma_m x) \le \varphi_{m+1}(x) + \gamma_m, \ x \ge 0.
$$

Then each function $f \in {\mathcal E}(\Psi^*)$ admits (the unique) extension to entire function belonging to $E(\varPhi)$.
\end{theorem}

The proofs of these theorems follows the standard schemes. 
Also they allow to obtain an additional information on a structure of the space $E(\varPhi)$. 
Namely, for each $\nu \in {\mathbb N}$ let 
${\mathcal H}(\varphi_{\nu})$ be a projective limit of spaces
$$
{\mathcal H}_k(\varphi_{\nu})= \{f \in H({\mathbb C}^n): {\cal N}_{\nu, k}(f) = 
\sup_{z \in {\mathbb C}^n} 
\frac 
{\vert f(z)\vert (1 + \Vert z \Vert)^k}
{e^{(\psi_{\nu}^*)^*(\ln (1 + \Vert Im z \Vert))}} < \infty \}, \ k \in {\mathbb Z}_+.
$$
Define the space ${\mathcal H}(\varPhi)$ as an inductive limit of spaces ${\mathcal H}(\varphi_{\nu})$.
In section 3 we show that if the family $\varPhi$ satisfies the condition
$i_3)$ then $E(\varPhi) = {\mathcal H}(\varPhi)$ (see Proposition 2).

For each $\nu \in {\mathbb N}$ and $m \in {\mathbb Z}_+$ define the normed space
$$
G_m({\psi_{\nu}^*}) = \{f \in C^m({\mathbb R}^n): 
\Vert f \Vert_{m, \psi_{\nu}^*}
 = \sup_{x \in {\mathbb R}^n, \vert \alpha \vert \le m, \beta \in {\mathbb Z}_+^n}  
\frac 
{\vert x^{\beta}(D^{\alpha}f)(x) \vert}
{\vert \beta \vert! e^{-\psi_{\nu}^*(\vert \beta \vert)}} < \infty \}.
$$
Let $G({\psi_{\nu}^*})= \bigcap \limits_{m = 0}^{\infty} G_m({\psi_{\nu}^*})$,
$G({\Psi^*})= \bigcup \limits_{\nu = 1}^{\infty} G({\psi_{\nu}^*})$. 
With usual operations of addition and multiplication by complex numbers $G({\psi_{\nu}^*})$ and $G({\Psi^*})$ are linear spaces. 
Endow $G({\psi_{\nu}^*})$ with a topology defined by the family of norms 
$\Vert \cdot \Vert_{m, \psi_{\nu}^*}$ ($m \in {\mathbb Z}_+$).
Supply $G({\Psi^*})$ with a topology of an inductive limit of spaces $G({\psi_{\nu}^*})$.

For $f \in E(\varPhi)$ define its Fourier transform $\hat f$ by the formula
$$
\hat f(x) = \int_{{\mathbb R}^n} f(\xi) e^{-i \langle x, \xi \rangle} \ d \xi , \ x \in {\mathbb R}^n.
$$

In section 4 we describe the space $E(\varPhi)$ in terms of Fourier transformation under additional conditions on $\varPhi$. 
 
\begin{theorem} 
Let functions of the family $\varPhi$ satisfy the following conditions: 

$i_3)$. for each $m \in {\mathbb N}$ there is a constant $a_m > 0$ such that 
$$
\varphi_m(2x) \le \varphi_{m+1}(x) + a_m, \ x \ge 0;
$$

$i_5)$. for each $m \in {\mathbb N}$ there exist numbers 
$h_m > 1$ and $l_m > 0$ such that 
$$
2 \varphi_m(x) \le \varphi_{m + 1}(h_k x) + l_m, \ x \ge 0.
$$

Then Fourier transformation ${\mathcal F}: f \in E(\varPhi) \to  \hat f$ establishes an isomorphism of spaces $E(\varPhi)$ and $G({\Psi^*})$.
\end{theorem}

Further, let ${\varPhi^*}=\{{\varphi_{\nu}^*}\}_{\nu=1}^{\infty}$. 
For each $\nu \in {\mathbb N}$ and $m \in {\mathbb Z}_+$ define the space
$$
GS_m(\varphi_{\nu}^*) = \{f \in C^m({\mathbb R}^n): q_{m, \nu}(f) = 
\sup_{x \in {\mathbb R}^n, \vert \alpha \vert \le m} 
\frac {\vert D^{\alpha}f)(x) \vert} {e^{-\varphi_{\nu}^*(\Vert x \Vert)}} < \infty \}.
$$
For each $\nu \in {\mathbb N}$ let
$GS({\varphi_{\nu}^*}) = \bigcap \limits_{m \in {\mathbb Z_+}}GS_m({\varphi_{\nu}^*})$. 
Let 
$GS({\varPhi^*})= \bigcup \limits_{\nu \in {\mathbb N}}GS({\varphi_{\nu}^*})$. 
Note that for each $\nu \in {\mathbb N}$ and $m \in {\mathbb Z}_+$
$$
q_{m, \nu}(f) \le q_{m + 1, \nu}(f), \ f \in GS_{m+1}(\varphi_{\nu}^*);  
$$
$$
q_{m, \nu + 1}(f) \le e^{C(\nu, 1)} q_{m, \nu}(f) , \ f \in GS_m(\varphi_{\nu}^*)).  
$$
Hence, $GS_{m+1}(\varphi_{\nu}^*)$ is continuously embedded in $GS_m(\varphi_{\nu}^*)$. 
Endow $GS({\varphi_{\nu}^*})$ with a topology defined by the family of norms $q_{\nu, m}$ ($m \in {\mathbb Z}_+$). 
From the last inequality it follows that $GS({\varphi_{\nu}^*})$ 
is continuously embedded in $GS({\varphi_{\nu + 1}^*})$. 
Supply $GS({\varPhi^*})$ with an inductive limit topology of spaces $GS({\varphi_{\nu}^*})$. 

The main result of the section 5 is the following theorem.

\begin{theorem}
Let functions of the family $\varPhi$ be convex and satisfy the condition $i_3)$ of Theorem 3. Then $G({\Psi^*}) =GS({\varPhi^*})$.
\end{theorem}

Note that if functions $\varphi_m$ are defined on $[0, \infty)$ by the formula 
$\varphi_m(x) = \Omega(2^m x)$ then the family $\varPhi$ satisfies the assumptions of Theorem 4.

\section{Auxiliary results}

In the proofs of Theorems results of the following three lemmas will be used.



\begin{lemma} 
Let $g$ be a real-valued continuous function on $[0, \infty)$ such that for some positive $a$ and $b$
$$
g(x) \ge a x - b, \ x \ge 0.
$$

Then 
$$
(g[e])^*(x) \le x \ln\frac {x}{a} - x + b, \ x > 0.
$$ 
\end{lemma}

{\bf Proof}. 
For all $x > 0$ we have
$$
(g[e])^*(x) =\sup_{y >0}(xy - g[e](y)) \le \sup_{y >0}(xy - a e^y) + b \le 
$$
$$
\le \sup_{y \in {\mathbb R}}(xy - a e^y) + b = 
x \ln\frac {x}{a} - x + b.
$$ 

\begin{corollary}
If $g \in {\cal B}$
then for each $M>0$ there exists a number $A_M>0$ such that for all $x > 0$
$$
(g[e])^*(x) \le x \ln\frac {x}{M} - x + A_M.
$$ 
\end{corollary}


{\bf Remark 1}. Note that if $g \in {\cal B}$ then by Corollary 1 for each $b>0$ the series \ $\displaystyle \sum_{j=0}^{\infty} \frac {e^{(g[e])^*(j)}}{b^j j!}$ and  
$
\displaystyle \sum_{\alpha \in {\mathbb Z_+^n}} 
\frac 
{e^{(g[e])^*(\vert \alpha \vert)}}
{b^{\vert \alpha \vert} \vert \alpha \vert!} 
$
are converging.

\begin{lemma} 
Let $u, v \in {\mathcal B}$ and there exist numbers $\tau>0$ and $C > 0$ such that
$$
2 u(x) \le v(x+\tau) + C, \ x \ge 0.
$$

Then there exists a number $A>0$ such that 
$$
v^*(x+y) \le u^*(x) + u^*(y) + \tau (x+y) + A, \ x, y \ge 0.
$$
\end{lemma} 

{\bf Proof}. 
For all $x, y, t \in [0, \infty)$ we have
$$
u^*(x) + u^*(y) \ge (x+y) t - 2 u(t).
$$
Using the assumption of lemma  we get
$$
u^*(x) + u^*(y) \ge (x+y) (t+\tau) - v(t+\tau) - C - \tau (x+y).
$$
Hence, 
\begin{equation}
u^*(x) + u^*(y) \ge \sup_{\xi \ge \tau} ((x+y)\xi - v(\xi)) - C - \tau (x+y).
\end{equation}
Further, 
$$
\displaystyle\sup_{0 \le \xi < \tau} ((x+y) \xi - v(\xi)) \le (x+y) \tau  - \inf \limits_{0 \le \xi < \tau} v(\xi) \le 
(x+y) \tau - \inf \limits_{\xi \ge 0} v(\xi).
$$
Since for each  $g \in {\mathcal B}$ we have that 
$g^*(t) \ge -\inf \limits_{\xi \ge 0} g(\xi)$ for all $t \ge 0$ then
$$
\displaystyle\sup_{0 \le \xi < \tau} ((x+y) \xi - v(\xi)) \le 
(x+y) \tau + u^*(x) + u^*(y)  + 2 \inf \limits_{\xi \ge 0} u(\xi) - \inf \limits_{\xi \ge 0} v(\xi).
$$ 
From this and inequality (1) we obtain that for all $x, y \ge 0$
$$
v^*(x+y) \le u^*(x) + u^*(y) + \tau (x+y) + A, 
$$
where $A = \max (C, 2 \inf \limits_{\xi \ge 0} u(\xi) - \inf \limits_{\xi \ge 0} v(\xi))$. 
The proof is complete.

\begin{lemma} 
Let real-valued continuous functions $u, v \in C[0, \infty)$ satisfy the following conditions: 

1. $\displaystyle \lim_{x \to + \infty} \frac {u[e](x)}{x}=
\displaystyle \lim_{x \to + \infty} \frac {v[e](x)}{x} = + \infty$;

2. there are constants $\sigma > 1$ and $\gamma > 0$ such that
$$
u(\sigma x) \le v(x) + \gamma, \ x \ge 0.
$$

Then  
$$
(u[e])^*(x) - (v[e])^*(x) \ge  x \ln \sigma - \gamma, \ x \ge 0.
$$ 
\end{lemma}

{\bf Proof}. Obviously, 
$
u[e](t+\ln \sigma) \le v[e](t) + \gamma, \ t \ge 0.
$
Then for $x \ge 0$
$$
(u[e])^*(x) - (v[e])^*(x) = \displaystyle \sup \limits_{t \ge 0}(xt - u[e](t)) - 
\displaystyle \sup \limits_{t \ge 0}(xt - v[e](t)) \ge 
$$
$$\ge 
\displaystyle \sup \limits_{t \ge 0}(xt - u[e](t)) - 
\displaystyle \sup \limits_{t \ge 0}(xt - u[e](t + \ln \sigma)) - \gamma = 
$$
$$=
\displaystyle \sup \limits_{t \ge 0}(xt - u[e](t)) - 
\displaystyle \sup \limits_{t \ge 0}(x(t + \ln \sigma) - u[e](t + \ln \sigma)) + x \ln \sigma - \gamma 
\ge  x \ln \sigma - \gamma.
$$ 

\section{Equivalent description of the space $E(\varPhi)$ \\
under additional conditions on $\varPhi$}

{\bf Proof of Theorem 1}. 
Let $f \in E(\varPhi)$. 
Then $f \in E(\varphi_{\nu})$ for some $\nu \in {\mathbb N}$. 
Let $m \in {\mathbb Z_+}$, $\alpha \in {\mathbb Z}_+^n$ and $x \in {\mathbb R}^n$ be arbitrary. 
Using Cauchy integral formula we have 
$$
(1+ \Vert x \Vert)^m (D^{\alpha}f)(x) =
\frac {\alpha! }{(2\pi i)^n} 
\displaystyle 
\idotsint \limits_{L_R(x)}
\frac 
{f(\zeta) (1+ \Vert x \Vert)^m \ d \zeta}
{(\zeta_1 - x_1)^{\alpha_1 +1} \cdots (\zeta_n - x_n)^{\alpha_n +1}} ,
$$
where for any
$R>0$ \ 
$L_R(x)= \{\zeta = (\zeta_1, \ldots , \zeta_n) \in {\mathbb C}^n: \vert \zeta_j - x_j \vert = R, j=1, \ldots , n \}$.
From this we get 
$$
(1+ \Vert x \Vert)^m \vert (D^{\alpha}f)(x) \vert \le 
$$
$$
\le 
\frac {\alpha! }{(2\pi)^n} 
\displaystyle \idotsint \limits_{L_R(x)}
\frac 
{(1+ \Vert x  - \zeta\Vert)^m (1+ \Vert \zeta \Vert)^m   \vert f(\zeta) \ \vert d \zeta \vert}
{\vert \zeta_1 - x_1\vert^{\alpha_1 +1} \cdots \vert \zeta_n - x_n\vert^{\alpha_n +1}}  \le
$$
$$
\le \frac 
{\alpha! p_{\nu, m}(f) (1 + nR)^m  e^{\varphi_{\nu}(nR)}}{R^{\vert \alpha \vert}} .
$$
Using the second condition on $\varPhi$ we have 
$$
(1+ \Vert x \Vert)^m \vert (D^{\alpha}f)(x) \vert \le 
e^{C(\nu, m)} \alpha! p_{\nu, m}(f)
\frac 
{e^{\varphi_{\nu + 1}(nR)}}{R^{\vert \alpha \vert}}.
$$
Hence, 
$$
(1+ \Vert x \Vert)^m \vert (D^{\alpha}f)(x) \vert \le 
e^{C(\nu, m)} n^{\vert \alpha \vert} \alpha! p_{\nu, m}(f)
\inf_{R>0}
\frac 
{e^{\varphi_{\nu + 1}(R)}}{R^{\vert \alpha \vert}}  \le 
$$
$$ 
\le e^{C(\nu, m)} n^{\vert \alpha \vert} \alpha! p_{\nu, m}(f)  
\exp(
{-\sup_{R>1} (\vert \alpha \vert \ln R - \varphi_{\nu + 1}(R)}) = 
$$
$$
= e^{C(\nu, m)} n^{\vert \alpha \vert} \alpha! p_{\nu, m}(f)
\exp({-\sup_{r > 0}
(\vert \alpha \vert r  - \psi_{\nu + 1} (r))}) = 
$$
$$
= e^{C(\nu, m)} n^{\vert \alpha \vert} \alpha! p_{\nu, m}(f) e^{-\psi_{\nu + 1}^*(\vert \alpha \vert)}.
$$
Note that by the condition $i_3)$ and Lemma 3 for each $k \in {\mathbb N}$
\begin{equation}
\psi_k^*(x) - \psi_{k + 1}^*(x) \ge x \ln 2 - a_k, \ x \ge 0,
\end{equation}
So using the inequality (2) we obtain that 
$$
(1+ \Vert x \Vert)^m \vert (D^{\alpha}f)(x) \vert \le 
a_{\nu, m} p_{\nu, m}(f) \alpha! e^{-\psi_{\nu + 2n +1}^*(\vert \alpha \vert)}, \ x \in {\mathbb R}^n, \alpha \in {\mathbb Z}_+^n,
$$
where $a_{\nu, m}$ is some positive number depending on $\nu$ and $m$.
Hence, we have that for each $m \in {\mathbb Z}_+$  
\begin{equation}
{\cal R}_{m, \nu + 2n + 1}(f_{|{\mathbb R}^n}) \le a_{\nu, m} p_{\nu, m}(f). 
\end{equation}
Therefore, 
$f_{|{\mathbb R}^n} \in {\mathcal E}(\psi_{\nu+2n+1}^*)$.
Thus, $f_{|{\mathbb R}^n} \in {\mathcal E}(\Psi^*)$.

{\bf Proof of Theorem 2}. 
Let $f \in {\mathcal E}(\Psi^*)$. 
Then $f \in {\mathcal E}(\psi_{\nu}^*)$ for some $\nu \in {\mathbb N}$.
Hence, for each $m \in {\mathbb Z}_+$ we have
\begin{equation}
(1+ \Vert x \Vert)^m \vert (D^{\alpha}f)(x) \vert \le  {\cal R}_{m, \nu} (f)
\alpha! e^{-\psi_{\nu}^*(\vert \alpha \vert)}, \ x \in {\mathbb R}^n, \alpha \in {\mathbb Z}_+^n.
\end{equation}
Since
$\displaystyle \lim_{x \to + \infty} \frac {\psi_{\nu}^*(x)}{x}= + \infty$ 
then for each $\varepsilon >0$ there is a constant $c_{\varepsilon}>0$ such that 
for all $x \in {\mathbb R}^n$ and $\alpha \in {\mathbb Z}_+^n$ \
$
\vert (D^{\alpha}f)(x)\vert \le c_{\varepsilon} {\varepsilon}^{\vert \alpha \vert}\alpha!. 
$
So it is clear that the sequence
$(\sum \limits_{\vert \alpha \vert \le k} \frac 
{(D^{\alpha}f)(0)}{\alpha!} x^{\alpha})_{k=1}^{\infty}$ 
converges to $f$ informly on compacts of ${\mathbb R}^n$ and the series 
$
\displaystyle \sum_{\vert \alpha \vert \ge 0} \frac {(D^{\alpha}f)(0)}{\alpha!} z^{\alpha}
$
converges informly on compacts of ${\mathbb C}^{n}$ and, hence, 
its sum $F_f(z)$ is an entire function in ${\mathbb C}^n$. 
Note that $F_f(x) = f(x), \ x \in {\mathbb R}^n$. The uniqueness of holomorphic continuation is obvious. 

Let us show that $F_f \in E(\varPhi)$. 
We will estimate a growth  of $F_f$ using the inequality (4) and the Taylor series expansion 
of $F_f$ with respect to a point $x \in {\mathbb R}^n$: 
$$
F_f(z) = \displaystyle \sum_{\vert \alpha \vert \ge 0} \frac {(D^{\alpha}f)(x)}{\alpha!} (iy)^{\alpha}, \ z = x+iy, 
y \in {\mathbb R}^n.
$$
Let $m \in {\mathbb Z}_+$ be arbitrary. Then
$$
(1 + \Vert z \Vert)^m \vert F_f(z) \vert \le \sum_{\vert \alpha \vert \ge 0} \frac 
{(1 + \Vert x \Vert)^m (1 + \Vert y \Vert)^{m + \vert \alpha \vert} \vert (D^{\alpha}f)(x)\vert }{\alpha!} 
\le 
$$
$$
\le 
\sum_{\vert \alpha \vert \ge 0}
{\cal R}_{m, \nu} (f) e^{-\psi_{\nu}^*(\vert \alpha \vert)}(1 + \Vert y \Vert)^{m + \vert \alpha \vert}  \le 
$$
$$
\le
{\cal R}_{m, \nu} (f) 
(1 + \Vert y \Vert)^m
\sum_{\vert \alpha \vert \ge 0} 
\frac {(1 + \Vert y \Vert)^{\vert \alpha \vert}}{e^{\psi^*_{\nu+1}(\vert \alpha \vert)}}e^{\psi^*_{\nu+1}(\vert \alpha \vert)-\psi^*_{\nu}(\vert \alpha \vert)}.
$$
Recall that $\varPhi$ satisfies the condition $i_4)$. 
Hence by Lemma 3 
\begin{equation}
\psi_k^*(x) - \psi_{k + 1}^*(x) \ge \delta_k x - \gamma_k, \ x \ge 0,
\end{equation}
where $\delta_k = \ln \sigma_k$.
Using this inequality and denoting $(\frac {e^{\gamma_{\nu} + \delta_{\nu}}}{e^{\delta_{\nu}} - 1})^n$ 
by $B_{\nu}$ we have  
$$
(1 + \Vert z \Vert)^m \vert F_f(z) \vert \le B_{\nu} {\cal R}_{m, \nu} (f) (1 + \Vert y \Vert)^m 
\sup \limits_{\vert \alpha \vert \ge 0}
\frac {(1 + \Vert y \Vert)^{\vert \alpha \vert}}{e^{\psi^*_{\nu+1}(\vert \alpha \vert)}} \le 
$$
$$
\le B_{\nu} {\cal R}_{m, \nu} (f) (1 + \Vert y \Vert)^m  e^{\sup \limits_{t \ge 0}(t \ln (1 + \Vert y \Vert) -  \psi^*_{\nu+1}(t))}.
$$ 
Thus,
\begin{equation}
(1 + \Vert z \Vert)^m \vert F_f(z) \vert 
\le B_{\nu} {\cal R}_{m, \nu} (f) 
e^{(\psi^*_{\nu+1})^*(\ln (1 + \Vert y \Vert)) + m \ln (1 + \Vert y \Vert))}.
\end{equation}

Note that from the condition $i_2)$ on $\varPhi$ it follows that 
for each $k \in {\mathbb N}$ and $A > 0$ 
$$
\psi_k(x) + A x \le \psi_{k+1}(x) + C(k, A), \ x \ge 0.
$$
So for all $\xi \ge 0$ 
$$
\psi_k^*(\xi) = \displaystyle \sup \limits_{x \ge 0}(\xi x - \psi_k(x)) \ge 
\displaystyle \sup \limits_{x \ge 0}(\xi x - \psi_{k + 1}(x) + A x) - C(k, A) = 
$$
$$
= \displaystyle \sup \limits_{x \ge 0}((\xi + A) x - \psi_{k + 1}(x)) - C(k, A) = 
\psi_{k + 1}^*(\xi + A) - C(k, A).
$$
Then for all $x \ge 0$ 
$$
(\psi_k^*)^*(x) =  \displaystyle \sup \limits_{x \ge 0}(x \xi  - \psi_k^*(\xi)) \le 
\displaystyle \sup \limits_{x \ge 0}(x \xi  - \psi_{k + 1}^*(\xi + A)) + C(k, A) = 
$$
$$
= \displaystyle \sup \limits_{x \ge 0}(x (\xi + A)  - \psi_{k + 1}^*(\xi + A)) - A x + C(k, A) \le 
(\psi_{k + 1}^*)^*(x) - A x + C(k, A).
$$
Thus, for each $k \in {\mathbb N}$ and $A > 0$ we have
\begin{equation}
(\psi_k^*)^*(x) + A x \le (\psi_{k+1}^*)^*(x) + C(k, A), \ x \ge 0.
\end{equation}
Now using the inequality  (7) we obtain from the estimate (6) that
\begin{equation}
(1 + \Vert z \Vert)^m \vert F_f(z) \vert 
\le B_{\nu} {\cal R}_{m, \nu} (f)  e^{C_{\nu + 1, m}} e^{(\psi_{\nu + 2}^*)^* (\ln (1 + \Vert y \Vert))}.
\end{equation}
It is clear that
$$
(1 + \Vert z \Vert)^m \vert F_f(z) \vert 
\le B_{\nu} {\cal R}_{m, \nu} (f)  e^{C_{\nu + 1, m}} e^{\psi_{\nu + 2}(\ln (1 + \Vert y \Vert))}.
$$
Hence,
$$
(1 + \Vert z \Vert)^m \vert F_f(z) \vert 
\le B_{\nu}
{\cal R}_{m, \nu} (f)  e^{C_{\nu + 1, m}} e^{\varphi_{\nu + 2} (1 + \Vert y \Vert)}.
$$
Using nondecreasity of functions of the family $\varPhi$ щт уфср мфкшфиду and the condition $i_2)$ on $\varPhi$ 
it is possible to find a constant $K_{\nu, m} > 0$ such that for all $z \in {\mathbb C}^n$
\begin{equation}
(1 + \Vert z \Vert)^m \vert F_f(z) \vert \le K_{\nu, m}
{\cal R}_{m, \nu} (f)  e^{\varphi_{\nu+3} (\Vert Im z \Vert)}.
\end{equation} 
Thus, for each $m \in {\mathbb Z}_+$ \ $p_{\nu + 3, m} (F_f) \le K_{\nu, m} {\cal R}_{m, \nu} (f) $. 
Hence, 
$F_f \in E(\varphi_{\nu + 3})$. Thus, $F_f \in E(\varPhi)$. 

Theorem 2 is proved.

Recall that in Section 1 we introduced the spaces ${\mathcal H}_k(\varphi_{\nu})$, 
${\mathcal H}(\varphi_{\nu})$ and ${\mathcal H}(\varPhi)$ as follows.
For each $\nu \in {\mathbb N}$ and $k \in {\mathbb Z}_+$ let
$$
{\mathcal H}_k(\varphi_{\nu})= \{f \in H({\mathbb C}^n): {\cal N}_{\nu, k}(f) = 
\sup_{z \in {\mathbb C}^n} 
\frac 
{\vert f(z)\vert (1 + \Vert z \Vert)^k}
{e^{(\psi_{\nu}^*)^*(\ln (1 + \Vert Im z \Vert))}} < \infty \}.
$$
Let ${\mathcal H}(\varphi_{\nu})= \bigcap \limits_{k=0}^{\infty} {\mathcal H}_k(\varphi_{\nu})$, 
${\mathcal H}(\varPhi)= \bigcup \limits_{\nu=1}^{\infty} {\mathcal H}(\varphi_{\nu})$. 
With usual operations of addition and multiplication by complex numbers 
${\mathcal H}(\varphi_{\nu})$ and ${\mathcal H}(\varPhi)$ are linear spaces. 
Since ${\cal N}_{\nu, k}(f) \le {\cal N}_{\nu, k+1}(f)$ for $f \in {\mathcal H}_{k+1}(\varphi_{\nu})$ 
then ${\mathcal H}_{k+1}(\varphi_{\nu})$ is continuously embedded in ${\mathcal H}_k(\varphi_{\nu})$. 
Endow ${\mathcal H}(\varphi_{\nu})$ with a projective limit topology of spaces ${\mathcal H}_k(\varphi_{\nu})$. 
Taking into account the inequality (7) we see that if $f \in {\mathcal H}(\varphi_{\nu})$ then 
${\cal N}_{\nu+1, k}(f) \le e^{C(\nu, 1)} {\cal N}_{\nu, k}(f)$ for each $k \in {\mathbb Z}_+$. 
Thus, ${\mathcal H}(\varphi_{\nu})$ is continuously embedded in ${\mathcal H}(\varphi_{\nu + 1})$ for each $\nu \in {\mathbb N}$. 
Supply ${\mathcal H}(\varPhi)$ with a topology of an inductive limit of spaces ${\mathcal H}(\varphi_{\nu})$.

\begin{proposition}
Let the family $\varPhi$ satisfies the condition $i_3)$. 
Then $E(\varPhi) = {\mathcal H}(\varPhi)$.
\end{proposition}

{\bf Proof}. 
First show that ${\mathcal H}(\varPhi)$ is continuously embedded in $E(\varPhi)$. 
Let $\nu \in {\mathbb N}$ be arbitrary and $f \in {\mathcal H}(\varphi_{\nu})$. 
Using nondecreasity of $\varphi_{\nu}$ and the condition $i_3)$ on $\varPhi$ 
we can find a constant $K_{\nu} > 0$ such that for each $k \in {\mathbb Z}_+$ 
$$
p_{\nu + 1, k}(f) \le K_{\nu} {\cal N}_{\nu, k}(f), \ f \in {\mathcal H}(\varphi_{\nu}).
$$
From this it follows that the embedding $I: {\mathcal H}(\varPhi) \to E(\varPhi)$ is continuous. 

The mapping $I$ is surjective too. Indeed, if $f \in E(\varPhi)$ then $f \in E(\varphi_{\nu})$ for some $\nu \in {\mathbb N}$. 
Let $m \in {\mathbb Z}_+$ be arbitrary. Recall that by the inequality (3)  \  
$
{\cal R}_{m, \nu + 2n + 1}(f_{|{\mathbb R}}) \le a_{\nu, m}  p_{\nu, m}(f).
$
From this and the inequality (8) (with $\nu$ replaced by $\nu + 2n + 1$; also recall that in our case $\sigma =2$ for each $m \in {\mathbb N}$) we obtain 
$$
{\cal N}_{\nu + 2n + 3, m}(f) \le A_{\nu, m} p_{\nu, m}(f),
$$
where $A_{\nu, m}$ is some positive number. Hence, $f \in {\mathcal H}(\varphi_{\nu + 2n + 3})$. 
So, $f \in {\mathcal H}(\varPhi)$. Moreover, the last estimate shows that the inverse mapping $I^{-1}$ is continuous.
Hence, the equality $E(\varPhi) = {\mathcal H}(\varPhi)$ is topological too.

\section{Fourier transformation of $E(\varPhi)$}

\hspace {0,5cm} First note that if $\varPhi$ satisfies the condition $i_3)$ then 
the space $G(\Psi^*)$ admits more simple description.  To show that we introduce the space $Q(\Psi^*)$ as follows.
For $\nu \in {\mathbb N}$ and $m \in {\mathbb Z}_+$ let 
$$
Q_m(\psi_{\nu}^*) = \{f \in C^m({\mathbb R}^n): 
$$
$$
N_{\nu, m}(f) = \max_{\vert \alpha \vert \le m}
\sup_{x \in {\mathbb R}^n, k  \in {\mathbb Z}_+}  
\frac 
{(1 + \Vert x \Vert)^k \vert (D^{\alpha}f)(x) \vert}
{k! e^{-\psi_{\nu}^*(k)}} < \infty \}.
$$
Let
$Q(\psi_{\nu}^*)= \bigcap \limits_{m \in {\mathbb Z_+}} Q_m(\psi_{\nu}^*)$ 
and 
$Q(\Psi^*)= \bigcup \limits_{\nu \in {\mathbb N}} Q(\psi_{\nu}^*).$
With usual operations of addition and multiplication by complex numbers $Q(\psi_{\nu}^*)$ and $Q(\Psi^*)$ are linear spaces. 
The family of norms $N_{\nu, m}(f)$ ($m \in {\mathbb Z}_+$) defines a locally convex topology in $Q(\psi_{\nu}^*)$. 
Endow $Q(\Psi^*)$ with a topology of inductive limit of spaces $Q(\psi_{\nu}^*)$. 

\begin{lemma}  
Let the family $\varPhi$ satisfies the condition $i_3)$. Then spaces $Q(\Psi^*)$ and $G(\Psi^*)$ coincide. 
\end{lemma} 
 
{\bf Proof}. Obviously, if $\nu \in {\mathbb N}$ and $f \in Q(\psi_{\nu}^*)$ 
then for each $m \in {\mathbb Z}_+$ we have 
$\Vert f \Vert_{m, \psi_{\nu}^*} \le N_{\nu, m}(f)$. Hence, $f \in G_m(\psi_{\nu}^*)$. Thus, if $f \in Q(\Psi^*)$ 
then $f \in G(\Psi^*)$ and the embedding mapping $J: Q(\Psi^*) \to G(\Psi^*)$ is continuous.

Show that $J$ is surjective. Indeed, let $\nu \in {\mathbb N}$, $f \in G(\psi_{\nu}^*)$ and 
$m \in {\mathbb N}$ be arbitrary.
Using the inequality (2)  we obtain
$$
\sup_{\Vert x \Vert \le 1, k  \in {\mathbb Z}_+}  
\frac {(1 + \Vert x \Vert)^k \vert (D^{\alpha}f)(x) \vert}
{k! e^{-\psi_{\nu + n}^*(k)}} 
\le 
\sup_{\Vert x \Vert \le 1, k  \in {\mathbb Z}_+}  
\frac {e^{a_{\nu+n-1}}\vert (D^{\alpha}f)(x) \vert}
{k! e^{-\psi_{\nu + n - 1}^*(k)}} .
$$
Take into account that $\lim \limits_{k \to \infty}\frac {e^{\psi_{\nu + n - 1}^*(k)}} {k!} = 0$ (see Corollary 1). 
So there is a number $C_1(\nu) > 1$ such that for each $\alpha \in {\mathbb Z_+^n}$
\begin{equation}
\sup_{\Vert x \Vert \le 1, k  \in {\mathbb Z}_+}  
\frac {(1 + \Vert x \Vert)^k \vert (D^{\alpha}f)(x) \vert}
{k! e^{-\psi_{\nu + n}^*(k)}} 
\le 
C_1(\nu) \sup_{\Vert x \Vert \le 1}  \vert (D^{\alpha}f)(x) \vert .
\end{equation}
Since for all $\alpha \in {\mathbb Z_+^n}$ with $\vert \alpha \vert \le m$ 
$$
\vert (D^{\alpha}f)(x) \vert \le \Vert f \Vert_{m, \psi_{\nu}^*} e^{-\psi_{\nu}^*(0)}, \ x \in {\mathbb R}^n,
$$
then from (10) we have that
\begin{equation}
\max_{\vert \alpha \vert \le m} \sup_{\Vert x \Vert \le 1, k  \in {\mathbb Z}_+}  
\frac 
{(1 + \Vert x \Vert)^k \vert (D^{\alpha}f)(x) \vert}{k! e^{-\psi_{\nu + n}^*(k)}} \le C_1(\nu)
\Vert f \Vert_{m, \psi_{\nu}^*} e^{-\psi_{\nu}^*(0)}.
\end{equation}
Further, for each $\alpha \in {\mathbb Z_+^n}$  
$$
\sup_{\Vert x \Vert > 1, k  \in {\mathbb Z}_+}  
\frac {(1 + \Vert x \Vert)^k \vert (D^{\alpha}f)(x) \vert}{k! e^{-\psi_{\nu + n}^*(k)}} \le 
\sup_{\Vert x \Vert > 1, \beta \in {\mathbb Z}_+^n}  
\frac 
{(2n)^{\vert \beta \vert}\vert x^{\beta}(D^{\alpha}f)(x) \vert}
{\vert \beta \vert! 
e^{-\psi_{\nu + n}^*(\vert \beta \vert)}}.
$$
Using the inequality (2) we get that 
for all $\alpha \in {\mathbb Z_+^n}$  
with $\vert \alpha \vert \le m$ that
$$
\sup_{\Vert x \Vert > 1, k  \in {\mathbb Z}_+}  
\frac {(1 + \Vert x \Vert)^k \vert (D^{\alpha}f)(x) \vert}{k! e^{-\psi_{\nu + n}^*(k)}}
\le 
\sup_{\Vert x \Vert > 1, \beta \in {\mathbb Z}_+^n}  
\frac 
{C_2(\nu) \vert x^{\beta}(D^{\alpha}f)(x) \vert}
{\vert \beta \vert! 
e^{-\psi_{\nu}^*(\vert \beta \vert)}} \le 
$$
$$
\le C_2(\nu)
\Vert f \Vert_{m, \psi_{\nu}^*}, 
$$
where $C_2(\nu) > 1$ is some constant.
From this and (11) we obtain for each $m \in {\mathbb Z}_+$ \ 
\begin{equation}
N_{\nu + n, m}(f) \le C(\nu) \Vert f \Vert_{m, \psi_{\nu}^*}, \ f \in G(\psi_{\nu}^*),
\end{equation}
where $C_1(\nu) = \max(C_2(\nu), C_1(\nu)e^{-\psi_{\nu}^*(0)})$. 
Hence, $f \in Q(\psi_{\nu + n}^*)$. Thus, if $f \in G(\Psi^*)$ 
then $f \in Q(\Psi^*)$. 
Note that from (12) it easily follows that the inverse mapping $J^{-1}$ is continuous. 
Therefore, the topological equality $Q(\psi^*)=G(\psi^*)$ is established.

{\bf Proof of Theorem 3}. Let $\nu \in {\mathbb N}$ and $f \in E(\varphi_{\nu})$.
Obviously, for all $\alpha$, $\beta \in {\mathbb Z}_+^n$, $x, \eta \in {\mathbb R}^n$
$$
x^{\beta} ({D^{\alpha} \hat f)(x) = x^{\beta} \int_{{\mathbb R}^n}} f(\zeta) 
(-i \zeta)^{\alpha} 
e^{-i \langle x, \zeta \rangle} \ d \xi, \  \zeta = \xi + i\eta.
$$
From this equality we have
$$
\vert x^{\beta} ({D^{\alpha} \hat f)(x)\vert \le 
\int_{{\mathbb R}^n}} 
\vert f(\zeta) \vert
\Vert \zeta \Vert^{\vert \alpha \vert} 
e^{\langle x, \eta \rangle} \Vert x \Vert^{\vert \beta \vert} \ d \xi \le 
$$
$$
\le 
\int_{{\mathbb R}^n} 
\vert f(\zeta) \vert
(1 + \Vert \zeta \Vert)^{n + \vert \alpha \vert + 1} 
e^{\langle x, \eta \rangle} \Vert x \Vert^{\vert \beta \vert} \ 
\frac {d \xi}{(1 + \Vert \xi \Vert)^{n+1}} \ .
$$
If $\vert \beta \vert = 0$ then 
\begin{equation}
\vert (D^{\alpha}\hat f)(x) \vert 
\le  s_n(1) e^{\varphi_{\nu} (0)} p_{\nu, n + \vert \alpha \vert + 1}(f).
\end{equation} 
If $\vert \beta \vert > 0$ and $x \ne 0$
then putting
$\eta = -\frac {x}{\Vert x \Vert} t$ \ ($t > 0$) we have
$$
\vert x^{\beta} 
(D^{\alpha}\hat f)(x)\vert \le s_n(1) p_{\nu, n + \vert \alpha \vert + 1}(f) 
e^{-t \Vert x \Vert} 
e^{\varphi_{\nu}(t)} \Vert x \Vert^{\vert \beta \vert}  \le 
$$
$$
\le s_n(1) p_{\nu, n + \vert \alpha \vert + 1}(f)
e^{\sup \limits_{r>0} (-tr + \vert \beta \vert \ln r)} e^{\varphi_{\nu}(t)}  =
$$
$$
= s_n(1) p_{\nu, n + \vert \alpha \vert + 1}(f)
e^{\vert \beta \vert \ln \vert \beta \vert - \vert \beta \vert - \vert \beta \vert \ln t} e^{\varphi_{\nu}(t)} .
$$
Since for each  $k \in {\mathbb Z}_+$
$$
\inf_{t>0} (- k \ln t + \varphi_{\nu}(t)) = - \sup_{t>0} (k \ln t - \varphi_{\nu}(t)) \le 
$$
$$
\le 
- \sup_{t \ge 1} (k \ln t - \varphi_{\nu}(t))= - \sup_{u \ge 0} (k u - \psi_{\nu}(u)) = - \psi_{\nu}^*(k),
$$
then from this and the previous estimate we obtain
\begin{equation}
\vert x^{\beta} 
(D^{\alpha} \hat f)(x)\vert \le s_n(1)
p_{\nu, n + \vert \alpha \vert + 1}(f)  
e^{\vert \beta \vert \ln \vert \beta \vert - \vert \beta \vert} 
e^{- \psi_{\nu}^*(\vert \beta \vert)}.
\end{equation}
If $\vert \beta \vert > 0$ and $x = 0$ then $x^{\beta} (D^{\alpha} \hat f)(x) =0$. 
From this and the inequalities (13) and (14) it follows that for all $\alpha, \beta \in {\mathbb Z}_+^n$, 
$x \in {\mathbb R}^n$
$$
\vert x^{\beta} 
(D^{\alpha} \hat f)(x)\vert \le s_n(1)
p_{\nu, n + \vert \alpha \vert + 1}(f) \vert \beta \vert! e^{- \psi_{\nu}^*(\vert \beta \vert)}.
$$
Thus, for each $m \in {\mathbb Z}_+$
$$
\max_{\vert \alpha \vert \le m} \sup_{x \in {\mathbb R}^n, \beta \in {\mathbb Z}^n_+}
\frac {\vert x^{\beta} 
(D^{\alpha} \hat f)(x)\vert}{\vert \beta \vert! e^{- \psi_{\nu}^*(\vert \beta \vert)}} \le
s_n(1) p_{\nu, n + m + 1}(f), \ f \in E(\varphi_{\nu}).
$$
In other words,
$$
\Vert \hat f \Vert_{m, \psi_{\nu}^*} \le s_n(1) p_{\nu, n + m + 1}(f), \ f \in E(\varphi_{\nu}).
$$
From this inequality it follows that the linear mapping ${\cal F}: f \in E(\varPhi)  \to \hat f$ acts from  
$E(\varPhi)$ to $G(\Psi^*)$ and is continuous. 

Let us show that ${\mathcal F}$ is surjective. 
Let $g \in G(\Psi^*)$. Then $g \in G(\psi_{\nu}^*)$ for some $\nu \in {\mathbb N}$. 
By the proof of Lemma 4 $g \in Q(\psi_{\nu + n}^*)$.
For all
$k \in {\mathbb Z}_+$, $\alpha \in {\mathbb Z_+^n}$  
$x \in {\mathbb R}^n$ we have
\begin{equation}
{(1 + \Vert x \Vert)^k \vert (D^{\alpha}g)(x) \vert} \le N_{\nu + n, \vert \alpha \vert}(g) 
k! e^{-\psi_{\nu + n}^*(k)}.
\end{equation} 
Let
$$
f(\xi) = 
\frac {1}{(2 \pi)^n} 
\int_{{\mathbb R}^n} g(x) e^{i \langle x, \xi \rangle} \ dx, \ 
\xi \in {\mathbb R}^n.
$$
For all $\alpha = (\alpha_1, \ldots , \alpha_n), \beta = (\beta_1, \ldots , \beta_n) \in {\mathbb Z_+^n}, \xi \in {\mathbb R}^n$ we have
$$
(i\xi)^{\beta}(D^{\alpha}f)(\xi) = \frac {(-1)^{\vert \beta \vert} }{(2 \pi)^n} 
\int_{{\mathbb R}^n} D^{\beta}(g(x) (ix)^{\alpha}) e^{i \langle x, \xi \rangle} \ dx .
$$
Put $\gamma_s =\min (\beta_s, \alpha_s)$ for $s=1, \ldots , n$, $\gamma = (\gamma_1, \ldots , \gamma_n)$.
Then 
$$
(i\xi)^{\beta}(D^{\alpha}f)(\xi) = \frac {(-1)^{\vert \beta \vert} }{(2 \pi)^n} 
\int_{{\mathbb R}^n}
\displaystyle \sum \limits_{j \in {\mathbb Z_+^n}: j \le \gamma}
C_{\beta}^j (D^{\beta - j} g)(x) (D^j (ix)^{\alpha}) 
e^{i \langle x, \xi \rangle} \ dx .
$$
From this (using the inequality (15)) we have
$$
\vert \xi^{\beta}(D^{\alpha}f)(\xi) \vert \le \frac {1}{(2 \pi)^n}
\displaystyle \sum \limits_{j \in {\mathbb Z_+^n}: j \le \gamma}
C_{\beta}^j 
\int_{{\mathbb R}^n} \vert(D^{\beta - j} g)(x)\vert 
\frac {\alpha!}{(\alpha - j)!} 
\Vert x \Vert^{\vert \alpha - j\vert}  \ dx \le
$$ 
$$
\le \frac {1}{(2 \pi)^n}
\displaystyle \sum \limits_{j \in {\mathbb Z_+^n}: j \le \gamma}
C_{\beta}^j \frac {\alpha!}{(\alpha - j)!}
\int \limits_{{\mathbb R}^n} \vert(D^{\beta - j} g)(x)\vert (1 + \Vert x \Vert)^{\vert \alpha  - j \vert + n+1} 
\frac {dx}{(1 + \Vert x \Vert)^{n+1}} \le
$$
$$
\le
\frac {s_n(1)}{(2 \pi)^n}
\displaystyle \sum \limits_{j \in {\mathbb Z_+^n}: j \le \gamma}
C_{\beta}^j
\frac {\alpha!}{(\alpha - j)!} 
\frac 
{N_{\nu + n, \vert \beta \vert}(g) (\vert \alpha \vert  - \vert j \vert + n+1)!} 
{e^{\psi_{\nu + n}^*(\vert \alpha \vert  - \vert j \vert + n+1)}} \le
$$
$$
\le
\frac {s_n(1) N_{\nu + n, \vert \beta \vert}(g) \alpha!}{(2 \pi)^n}
\displaystyle \sum \limits_{j \in {\mathbb Z_+^n}: j \le \gamma} 
C_{\beta}^j
\frac {(\vert \alpha \vert  - \vert j \vert + n+1)!}{(\alpha - j)!} 
e^{-\psi_{\nu + n}^*(\vert \alpha \vert  - \vert j \vert)} \ . 
$$
Note that from the condition $i_5)$ on $\varPhi$ it follows that for each $k \in {\mathbb N}$
$$
2\psi_k(x) \le \psi_{k+1}(x+b_k) + l_k, \ x \ge 0,
$$
where $b_k = \ln h_k$.
Hence, by Lemma 2 
\begin{equation}
\psi_{k+1}^*(x+y) \le \psi_k^*(x) + \psi_k^*(y) + b_k (x+y) + A_k, \ x, y \ge 0,
\end{equation}
where $A_k$ is some positive constant.
Using the inequality (16) and setting $c_1= \frac 
{s_n(1) e^{A_{\nu + n}}}{(2 \pi)^n}$ we have
$$
\vert \xi^{\beta}(D^{\alpha}f)(\xi) \vert \le 
\frac 
{c_1 e^{b_{\nu + n} \vert \alpha \vert} N_{\nu + n, \vert \beta \vert}(g) \alpha!} {e^{\psi_{\nu + n + 1}^*(\vert \alpha \vert)}}
\displaystyle \sum \limits_{j \in {\mathbb Z_+^n}: \atop j \le \gamma} 
\frac 
{C_{\beta}^j(\vert \alpha \vert  - \vert j \vert + n+1)!e^{\psi_{\nu + n}^*(\vert j \vert)}}{(\alpha - j)!} .
$$
Note that $
(m_1 + m_2)! \le e^{m_1 + m_2} m_1! m_2! 
$
for all $m_1, m_2 \in {\mathbb Z_+}$ 
From this it follows that for all
$m_1, \ldots , m_n \in {\mathbb Z_+}$ 
\begin{equation}
(m_1 + \cdots + m_n)! \le e^{(n-1)(m_1 + \cdots + m_n)} m_1! \cdots  m_n!. 
\end{equation}
Using this inequality and setting
$c_2 = c_1 e^{n+1}(n+1)!$ we have
$$
\vert \xi^{\beta}(D^{\alpha}f)(\xi) \vert \le 
\frac 
{c_2 e^{b_{\nu + n} \vert \alpha \vert} N_{\nu + n, \vert \beta \vert}(g)\alpha!} {e^{\psi_{\nu + n + 1}^*(\vert \alpha \vert)}}
\displaystyle \sum \limits_{j \in {\mathbb Z_+^n}: j \le \gamma} 
C_{\beta}^j
\frac 
{e^{\vert \alpha \vert  - \vert j \vert} (\vert \alpha \vert  - \vert j \vert)! 
e^{\psi_{\nu + n}^*(\vert j \vert)}}{(\alpha - j)!}.
$$
Using the inequality (17) again we obtain
$$
\vert \xi^{\beta}(D^{\alpha}f)(\xi) \vert 
\le 
\frac 
{c_2 e^{b_{\nu + n} \vert \alpha \vert} N_{\nu + n, \vert \beta \vert}(g)\alpha!} {e^{\psi_{\nu + n + 1}^*(\vert \alpha \vert)}}
\displaystyle \sum \limits_{j \in {\mathbb Z_+^n}: j \le \gamma} 
C_{\beta}^j
e^{n (\vert \alpha \vert  - \vert j \vert)} 
e^{\psi_{\nu + n}^*(\vert j \vert)}.
$$
Since the series 
$\displaystyle \sum_{\vert j \vert \ge 0} 
e^{- n \vert j \vert} \frac {e^{\psi_{\nu + n}^*(\vert j \vert)}}{j!}$ converges (see Remark 1) then we have
$$
\vert \xi^{\beta}(D^{\alpha}f)(\xi) \vert 
\le \frac 
{c_2 \beta!e^{(b_{\nu + n} + n) \vert \alpha \vert} N_{\nu + n, \vert \beta \vert}(g)\alpha!} 
{e^{\psi_{\nu + n + 1}^*(\vert \alpha \vert)}} 
\displaystyle \sum_{\vert j \vert \ge 0} 
e^{- n \vert j \vert} 
\frac 
{e^{\psi_{\nu + n}^*(\vert j \vert)}}
{j!}.
$$
Using (17) once more and setting
$c_3 = c_2 \beta! \displaystyle \sum_{\vert j \vert \ge 0} 
\frac 
{e^{\psi_{\nu + n}^*(\vert j \vert)}}
{\vert j \vert!}$ 
we have
$$
\vert \xi^{\beta}(D^{\alpha}f)(\xi) \vert 
\le 
c_3 e^{(b_{\nu + n} + n) \vert \alpha \vert} N_{\nu + n, \vert \beta \vert}(g)\alpha! {e^{-\psi_{\nu + n + 1}^*(\vert \alpha \vert)}} .
$$
From this using the inequality (2) we can find integer $s=s(\nu, n) > n + 1$ 
and a constant $c_4 > 0$ (depending on $\nu, n$ and $\beta$)
such that
$$
\vert \xi^{\beta}(D^{\alpha}f)(\xi) \vert 
\le 
c_4 N_{\nu + n, \vert \beta \vert}(g)\alpha! {e^{-\psi_{\nu + s}^*(\vert \alpha \vert)}} .
$$
So if $m \in {\mathbb Z}_+$ then from the last inequality we get that
for all $\alpha \in {\mathbb Z_+^n}, \xi \in {\mathbb R}^n$
$$
(1 + \Vert \xi \Vert)^m \vert (D^{\alpha}f)(\xi) \vert  \le 
c_5 N_{\nu + n, m}(g) \alpha! e^{-\psi_{\nu + s}^*(\vert \alpha \vert)},
$$
where $c_5 > 0$ is some positive constant depending on $\nu, n$ and $m$.
By Theorem 2 $f$ can be holomorphically continued (uniquely) to entire function $F_f$ belonging to $E(\varPhi)$.   
Obviously, 
$g = {\mathcal F}(F_f)$. 
The proof of Theorem 2 (see inequalities (4) and (9)) indicates that  there is a constant 
$c_6 > 0$  (depending on $\nu, n$ and $m$) such that for $z \in {\mathbb C}^n$
$$
(1 + \Vert z \Vert)^m \vert F_f(z) \vert \le c_6 N_{\nu + n, m}(g) e^{\varphi_{\nu + s + 3} (\Vert Im z \Vert)}.
$$
Hence, 
$$
p_{\nu + s + 3, m}(F_f) \le c_6 N_{\nu + n, m}(g). 
$$
Taking into account (12) we get 
$$
p_{\nu + s + 3, m}(F_f) \le c_7 \Vert g \Vert_{m, \psi_{\nu}^*}, \ g \in G(\psi_{\nu}^*),
$$
where $c_7 > 0$ depends on on $\nu, n$ and $m$.
From this estimate it follows
that the inverse mapping ${\mathcal F}^{-1}$ is continuous. 

Thus, we have proved that Fourier transformation establishes a topological isomorphism between spaces $E(\varPhi)$ and $G(\Psi^*)$.

\section{Special case of $\varPhi$}

In the proof of Theorem 4 the following three lemmas will be used.

\begin{lemma} 
Let $g$ be a real-valued continuous function on  
$[0, \infty)$ such that
$
\displaystyle \lim_{x \to + \infty} \frac {g(x)}{x}= + \infty.
$
Then for each $\delta > 0$
$$
\displaystyle \lim_{x \to + \infty} \frac {g^*((1+\delta) x) - g^*(x)}{x}= + \infty.
$$
\end{lemma} 

{\bf Proof}. Let $\delta > 0$ be arbitrary. For each $x>0$ denote by $\xi(x)$ a point where the supremum  of the function $u_x(\xi)= x \xi - g(\xi)$ over $[0, \infty)$ is attained. 
Note that  $\xi(x) \to +\infty$ as $x \to +\infty$. Otherwise there are a number $M>0$ and a sequence $(x_j)_{j=1}^{\infty}$ of positive numbers $x_j$ converging to $+\infty$ such that $\xi(x_j) \le M$. 
Then $g^*(x_j)= x_j \xi(x_j) - g(\xi(x_j))$.  But it contradicts to the fact that 
$
\displaystyle \lim_{x \to + \infty} \frac {g^*(x)}{x}= + \infty.
$
Thus, 
$
\displaystyle \lim_{x \to + \infty} \xi(x)= + \infty.
$
From this and the inequality
$$
g^*((1+\delta) x)- g^*(x) \ge (1+\delta) x \xi(x) - g(\xi(x)) - x \xi(x) + g(\xi(x))=\delta x \xi(x), \ x>0,
$$
the assertion of lemma follows. The proof is complete.

The following lemma follows from a result of S.V. Popenov \cite {N-P} (see Lemma 4 there). So the proof is omitted. 

\begin{lemma}
Let $u \in V$. Then there exists a constant $K > 0$ depending on $u$ such that 
$$
(u[e])^*(t) + (u^*[e])^*(t) \ge  t \ln t - t - K, \ t > 0.
$$
\end{lemma}

The following lemma was proved in \cite {M-P}, \cite {M}.

\begin{lemma}
Let a lower semi-continuous function $u: [0, \infty) \to {\mathbb R}$ be such that 
$\displaystyle \lim_{x \to + \infty} \frac {u(x)}{x}= + \infty$.
Then
$$
(u[e])^*(x) + (u^*[e])^*(x) \le x \ln x - x, \
x > 0.
$$
\end{lemma}

{\bf Proof of Theorem 4}.
Let $\nu \in {\mathbb N}$, $f \in G(\psi_{\nu}^*)$. Fix $m \in {\mathbb Z}_+$.
Since $f \in Q(\psi_{\nu + n}^*)$ (see the proof of Lemma 4) then 
for all $k \in {\mathbb Z}_+$, $\alpha \in {\mathbb Z_+^n}$ with $\vert \alpha \vert \le m$, 
$x \in {\mathbb R}^n$ 
\begin{equation}
\vert (D^{\alpha} f)(x) \vert \le 
\frac 
{N_{\nu + n, m}(f) 
k! e^{-\psi_{\nu + n}^*(k)}}
{(1 + \Vert x \Vert)^k} \ .
\end{equation}
Taking into account that $j! < \frac {3j^{j+1}}{e^j}$ for all $j \in {\mathbb N}$ 
and using the inequality (16) and nondecreasity of $\psi_{\nu + n}^*$ we have
for all $k \in {\mathbb N}$, $t \in [k, k+1)$ and $\mu \ge 1$
$$
\frac 
{k! e^{-\psi_{\nu + n}^*(k)}}{{\mu}^k} \le 
\frac {3 k^{k+1} e^{-\psi_{\nu + n}^*(k)}} {e^k {\mu}^k} \le 
\frac 
{3 \mu t^{t+1} e^{-\psi_{\nu + n + 1}^*(t) + \psi_{\nu + n}^*(1) + b_{\nu + n} t + A_{\nu + n}}}{e^t {\mu}^t}.
$$
Using the inequality (2) we have 
$$
\frac 
{k! e^{-\psi_{\nu + n}^*(k)}}{{\mu}^k} \le C_1 {\mu}
e^{(t+1) \ln t - \psi_{\nu + s}^*(t)  - t \ln e {\mu} + (b_{\nu + n} - s\ln 2) t},
$$
where integer $s \ge n + 2$, positive constant $C_1$ depends on $\nu, n$ and $s$.
Now choose $s \in {\mathbb N}$ so that $s\ln 2 > b_{\nu + n}$ 
then we can find a positive constant $C_2$ (depending on $\nu, n$ and chosen $s$) such that 
$$
\frac 
{k! e^{-\psi_{\nu + n}^*(k)}}{{\mu}^k} \le C_2 {\mu}
e^{t \ln t - \psi_{\nu + s}^*(t)  - t \ln {\mu}}.
$$
Now with help of Lemma 6 we obtain
$$
\frac 
{k! e^{-\psi_{\nu + n}^*(k)}}{{\mu}^k} \le C_3 {\mu}
e^{(\varphi_{\nu + s}^*[e])^*(t)  - t \ln {\mu}},
$$
where $C_3$ is some positive constant depending on $\nu, n$ and chosen $s$.
From this it follows that 
\begin{equation}
\displaystyle \inf \limits_{k \in {\mathbb N}} \frac 
{k! e^{-\psi_{\nu + n}^*(k)}}{{\mu}^k} \le 
C_3 \mu 
e^{\inf \limits_{t \ge 1}(\varphi_{\nu + s}^*[e])^*(t) - t \ln \mu}.
\end{equation}
Obviously, 
$$
\inf_{t \ge 1} \left((\varphi_{\nu + s}^*[e])^*(t)  - t \ln \mu \right) \le 
-\ln \mu + (\varphi_{\nu + s}^*[e])^*(1);
$$
$$
\inf_{0 < t \le 1} \left((\varphi_{\nu + s}^*[e])^*(t)  - t \ln \mu \right) \ge 
-\ln \mu + (\varphi_{\nu + s}^*[e])^*(0).
$$
Consequently,
$$
\inf_{t \ge 1} \left((\varphi_{\nu + s}^*[e])^*(t)  - t \ln \mu \right) \le 
$$
$$
\le
\inf_{0 < t \le 1} \left((\varphi_{\nu + s}^*[e])^*(t)  - t \ln \mu \right)  + (\varphi_{\nu + s}^*[e])^*(1) - (\varphi_{\nu + s}^*[e])^*(0).
$$
Denoting $(\varphi_{\nu + s}^*[e])^*(1) - (\varphi_{\nu + s}^*[e])^*(0)$ by $m_{\nu}$
we have 
$$
\inf_{t \ge 1} \left((\varphi_{\nu + s}^*[e])^*(t)  - t \ln \mu \right) \le 
\inf_{t > 0} \left((\varphi_{\nu + s}^*[e])^*(t)  - t \ln \mu \right)  + m_{\nu}.
$$
Going back to (19) we have
\begin{equation}
\inf_{k \in {\mathbb N}} \frac 
{k! e^{-\psi_{\nu + n}^*(k)}}{{\mu}^k} \le 
C_3 e^{m_{\nu}} {\mu}
e^{\inf \limits_{t > 0} \left((\varphi_{\nu + s}^*[e])^*(t)  - t \ln \mu \right)}.
\end{equation} 
For each $j \in {\mathbb N}$ choose $\theta_j \in V_{\varphi_j^*[e]}$.  
Then 
\begin{equation}
\vert \theta_{j}(\xi) - 
\varphi_j^*[e](\xi) \vert \le r_j, \ \xi \ge 0; 
\end{equation} 
\begin{equation}
\vert \theta_j^*(\xi) - (\varphi_j^*[e])^*(\xi) \vert \le r_j, \ \xi \ge 0,
\end{equation} 
where $r_j$ is some positive number depending on $\varphi_j^*[e]$ and $\theta_j$.
From (20) using the inequality (22) we have 
$$
\inf_{k \in {\mathbb N}} \frac 
{k! e^{-\psi_{\nu + n}^*(k)}}{{\mu}^k} \le C_4 \mu
e^{\inf \limits_{t > 0} \left(\theta_{\nu + s}^*(t)  - t \ln \mu \right)}.
$$
where
$C_4=C_3 e^{m_{\nu} + r_{\nu + s}}$.
Using the Young inversion formula we obtain
$$
\inf_{k \in {\mathbb N}} \frac 
{k! e^{-\psi_{\nu + n}^*(k)}}{{\mu}^k} \le C_4 \mu
e^{-\theta_{\nu + s} (\ln \mu)}.
$$
From this using the inequality (21) we have
$$
\inf_{k \in {\mathbb N}} \frac 
{k! e^{-\psi_{\nu + n}^*(k)}}{{\mu}^k} \le C_5 \mu 
e^{-\varphi_{\nu + s}^*[e] (\ln \mu)},  
$$
where
$C_5=C_{\nu, 4} e^{r_{\nu + s}}$.
In other words, 
$$
\inf_{k \in {\mathbb N}} \frac 
{k! e^{-\psi_{\nu + n}^*(k)}}{{\mu}^k} \le C_5 \mu 
e^{-\varphi_{\nu + s}^*(\mu)}. 
$$
Using this inequality and nondecreasity of $\varphi_{\nu + s}^*$ we have
\begin{equation}
\inf_{k \in {\mathbb N}} \frac 
{k! e^{-\psi_{\nu + n}^*(k)}}{(1 + \Vert x \Vert)^k} \le C_5
(1 + \vert x \vert) 
e^{-\varphi_{\nu + s}^*(\Vert x \Vert)}, \ x \in {\mathbb R}^n.
\end{equation}

Note that using the condition $i_3)$ on $\varPhi$ it is easy to obtain that for each $j \in {\mathbb N}$
\begin{equation}
\varphi_{j+1}^*(\xi) \le \varphi_j^*\left(\frac {\xi}{2}\right) + a_j, \ \xi \ge 0.
\end{equation}
Hence,
$
\varphi_j^*(\xi) - \varphi_{j+1}^*(\xi) \ge \varphi_j^*(\xi) - 
\varphi_j^*\left(\frac {\xi}{2}\right) - a_j, \ \xi \ge 0.
$
From this and Lemma 5 we get that
\begin{equation}
\displaystyle \lim_{x \to + \infty} \frac 
{\varphi_j^*(\xi) - \varphi_{j+1}^*(\xi)}{\xi}= + \infty.
\end{equation}

Going back to (23) we obtain with help of (25) that
$$
\inf_{k \in {\mathbb N}} \frac 
{k! e^{-\psi_{\nu +n}^*(k)}}{(1 + \Vert x \Vert)^k} \le C_6
e^{-\varphi_{\nu + s + 1}^*(\Vert x \Vert)}, \ 
x \in {\mathbb R}^n,
$$
where
$C_6$ is some positive number.

From this and the inequality (18) we obtain that for all 
$\alpha \in {\mathbb Z_+^n}$ with $\vert \alpha \vert \le m$
\begin{equation}
\vert (D^{\alpha} f)(x) \vert \le C_6 N_{\nu + n, m}(f)
e^{-\varphi_{\nu + s + 1}^*(\Vert x \Vert)}.
\end{equation}

This means that  
$$
q_{m, \nu + s + 1}(f) \le C_6 N_{\nu + n, m}(f), \ f \in G(\psi_{\nu}^*).
$$
Taking into account the inequality (12) we have
$$
q_{m, \nu + s + 1}(f) \le C_7 \Vert f \Vert_{m, \psi_{\nu}^*}, \ f \in G(\psi_{\nu}^*),
$$
where $C_7$ is some positive constant depending on $\nu$.
From this it follows that the identity mapping $I$ acts from $G(\Psi^*)$ to $GS(\varPhi^*)$ and is continuous.

Show that $I$ is surjective. Let $f \in GS(\varPhi^*)$. 
Then $f \in GS(\varphi_{\nu}^*)$ for some $\nu \in {\mathbb N}$. 
Let $m \in {\mathbb Z}_+$ be fixed, $x \in {\mathbb R}^n$ be arbitrary. 
For all $\alpha \in {\mathbb Z_+^n}$ such that $\vert \alpha \vert \le m$ we have 
\begin{equation}
\vert (D^{\alpha}f)(x) \vert \le q_{m, \nu}(f)
e^{-\varphi_{\nu}^*(\Vert x \Vert)}.
\end{equation}
Using the inequality (24) we have from (27) that 
$$
\vert (D^{\alpha}f)(x) \vert \le e^{a_{\nu}} q_{m, \nu}(f)
e^{-\varphi_{\nu + 1}^*(2 \vert x \vert)}.
$$
Obviously, there exists a constant $M_{\nu} > 1$ such that 
$$
\vert (D^{\alpha}f)(x) \vert   \le M_{\nu} q_{m, \nu}(f)
e^{-\varphi_{\nu + 1}^*(\Vert x \Vert + 1)}.
$$
In other words,
$$
\vert (D^{\alpha}f)(x) \vert   \le M_{\nu} q_{m, \nu}(f)
e^{-\varphi_{\nu + 1}^*[e](\ln (\Vert x \Vert + 1))}.
$$
From this we have
$$
\vert (D^{\alpha}f)(x) \vert  \le M_{\nu} q_{m, \nu}(f)  
e^{- \sup \limits_{t > 0} (t \ln (\Vert x \Vert + 1) - (\varphi_{\nu + 1}^*[e])^*(t))}.
$$
Now using Lemma 7 we get
$$
\vert (D^{\alpha}f)(x) \vert \le M_{\nu} q_{m, \nu}(f)  
e^{-\sup\limits_{t > 0} (t  \ln (e (\vert x \vert + 1))  - t \ln t + \psi_{\nu + 1}^*(t))}.
$$
Consequently, for all $k \in {\mathbb N}$
$$
\vert (D^{\alpha}f)(x) \vert \le M_{\nu} q_{m, \nu}(f)
\frac 
{k^k  e^{-\psi_{\nu + 1}^*(k)}}
{(e (1 + \Vert x \Vert))^k}.
$$
From this and (27) it follows that 
$$
(1 + \Vert x \Vert)^k
\vert (D^{\alpha}f)(x) \vert  \le M_{\nu}
q_{m, \nu}(f) 
k!  e^{-\psi_{\nu + 1}^*(k)}, \ k \in {\mathbb Z_+}.
$$
This means that 
\begin{equation}
\Vert f \Vert_{m, \psi_{\nu + 1}^*} \le M_{\nu} q_{m, \nu}(f).
\end{equation}
Since here $m \in {\mathbb Z}_+$ is arbitrary then $f \in G(\psi_{\nu + 1}^*)$. Hence, $f \in G(\Psi^*)$.
From (28) it follows that the mapping $I^{-1}$ is continuous. Thus, 
the spaces $G({\Psi^*})$ and $GS({\varPhi^*})$ coincides.

{\bf Acknowledgements}. 
The research of the first author was supported by grants from RFBR (14-01-00720, 14-01-97037).

\pagebreak


\begin{thebibliography}{99}

\bibitem {GS1} 
Gelfand I.M., Shilov G.E., {\it Generalized functions}, Vol. 2, Academic Press, New York,
1967.

\bibitem {GS2} 
Gelfand I.M., Shilov G.E., {\it Generalized functions}, Vol. 3, Academic Press, New York, 1967. 

\bibitem {C-C-K 1}  J. Chung, S.-Y. Chung, and D. Kim. 
{\it Characterizations of the Gelfand-Shilov spaces via Fourier transforms} // Proc. Amer. Math. Soc. 1996. V.~124. N. 7. P.~2101--2108.

\bibitem {C-C-K 2} J. Chung, S.-Y. Chung, and D. Kim. {\it Equivalence of the Gelfand-Shilov spaces} // Journal of Math. Anal. and Appl. 1996. V.~203. P.~828--839.

\bibitem {Gur1} 
Gurevich B.L., {\it New types of spaces of fundamental and generalized
functions and Cauchy's problem for systems of finite difference equations},
Doklady Akad. Nauk SSSR (N.S.), {\bf 99} (1954), 893--895.

\bibitem {Gur2}
Gurevich B.L., {\it New types of fundamental and generalized spaces and
Cauchy's problem for systems of difference equations involving differential
operations}, Doklady Akad. Nauk SSSR (N.S.), {\bf 108} (1956), 1001--1003.

\bibitem {M-P} 
Musin I.Kh., Popenov S.V., {\it On a weighted space of infinitely differentiable functions in ${\mathbb R}^n$}, Ufa Mathematical Journal, {\bf 2} (2010), 54--62. 

\bibitem {M} 
Musin I.Kh., {\it Approximation by polynomials in a weighted space of infinitely differentiable functions with an application to hypercyclicity}. Extracta Mathematicae, {\bf 27} (2012), 75--90. 

\bibitem{N-P}
Napalkov V.V., Popenov S.V.: On Laplace transformation on weighted Bergman space of entire functions 
on ${\mathbb C}^n$. Doklady Mathematics.  {\bf 55}, 110--112 (1997)

\bibitem {R-V}
Roberts A.W., Varberg D.E., {\it Convex functions}, Academic Press, New York and London, 1973.


\end{thebibliography}
\end{document}